\documentclass[reqno]{amsart}
\usepackage{amssymb,eucal,latexsym,enumerate,color,url,arydshln,mathdots,tikz,tikz-cd}

\newenvironment{enumeratei}{\begin{enumerate}[\upshape (i)]}%
{\end{enumerate}}
{\end{enumerate}}
\newenvironment{enumerater}{\begin{enumerate}[\upshape (1)]}%
{\end{enumerate}}

\newcommand{\pup}[1]{\textup{(}{#1}\textup{)}}

\newcommand{\eqdef}{\overset{\mathrm{def}}{=}}

\newcommand{\simr}{\sim_{\mathrm{r}}}
\newcommand{\siml}{\sim_{\mathrm{l}}}

\DeclareMathOperator{\End}{End}

\newcommand{\ga}{\alpha}
\newcommand{\gb}{\beta}

\newcommand{\gf}{\varphi}

\newcommand{\gl}{\lambda}

\newcommand{\gs}{\sigma}

\newcommand{\eps}{\varepsilon}

\newcommand{\gS}{\Sigma}

\newcommand{\ol}[1]{\overline{#1}}

\newcommand{\pI}[1]{\bigl({#1}\bigr)}

\newcommand{\set}[1]{\{#1\}}
\newcommand{\setm}[2]{\set{{#1}\mid{#2}}}
\newcommand{\vecm}[2]{({#1}\mid{#2})}
\newcommand{\seq}[1]{\langle{#1}\rangle}

\newcommand{\op}{\mathsf{op}}

\newcommand{\id}{\mathrm{id}}

\newcommand{\res}{\mathbin{\restriction}}

\newcommand{\FF}{\mathbb{F}}

\newcommand{\HH}{\mathbb{H}}

\newcommand{\ZZ}{\mathbb{Z}}

\newcommand{\RR}{\mathbb{R}}

\newcommand{\cS}{{\mathcal{S}}}

\numberwithin{equation}{section}

\newtheorem*{stat}{\name}
\newcommand{\name}{testing}

\theoremstyle{plain}

\newtheorem{theorem}{Theorem}[section]
\newtheorem{proposition}[theorem]{Proposition}
\newtheorem{corollary}[theorem]{Corollary}
\newtheorem{lemma}[theorem]{Lemma}
\newtheorem{examplepf}[theorem]{Example}
\newtheorem{claim}{Claim}
\newtheorem*{sclaim}{Claim}

\theoremstyle{definition}

\newtheorem{definition}[theorem]{Definition}

\newtheorem{example}[theorem]{Example}

\theoremstyle{remark}
\newtheorem{remark}[theorem]{Remark}
\newtheorem*{note}{Note}

\newcommand{\qedc}{{\qed}~{\rm Claim~{\theclaim}.}}
\newcommand{\qedsc}{{\qed}~{\rm Claim.}}

\newenvironment{cproof}
{\begin{proof}[Proof of Claim.]}
{\qedc\renewcommand{\qed}{}\end{proof}}

\newenvironment{scproof}
{\begin{proof}[Proof of Claim.]}
{\qedsc\renewcommand{\qed}{}\end{proof}}

\numberwithin{figure}{section}
\numberwithin{table}{section}

\newcommand{\xx}{\mathbf{x}}

\newcommand{\kk}{\Bbbk}
\newcommand{\kkm}{\kk^{\times}}
\newcommand{\FFm}{\FF^{\times}}

\newcommand{\ls}[2]{{#1}(\!({#2})\!)}

\author[F. Wehrung]{Friedrich Wehrung}
\address{Normandie Universit\'e, UNICAEN\\
CNRS UMR 6139, LMNO\\
14000 Caen\\
France}
\email{friedrich.wehrung01@unicaen.fr}
\urladdr{https://wehrungf.users.lmno.cnrs.fr}

\title[Similarity in domains]{A note on similarity classes from atomic factorizations in domains}

\date{\today}



\begin{document}

\begin{abstract}
We construct an example of a non-Archimedean valued skew field~$\kk$, with elements~$a$ and~$b$ of distinct absolute value such that there are no \linebreak$p,q\in\kk[\xx]$ ($\xx$ a central indeterminate) such that $(\xx-b)(\xx-a)=pq$, $p$ similar to~$\xx-a$, and~$q$ similar to~$\xx-b$.
The obstruction lies in the non-existence of a solution in~$\kk$ to Sylvester's equation $at-tb=1$.
We also disprove a 1963 result by Cohn (also reproduced in Cohn's 2006 monograph), by constructing a domain with similar elements~$a$ and~$b$ with~$a$ an atom and~$b$ not a product of atoms (resp., a product of two atoms, one of them central).
\end{abstract}

\maketitle

\tableofcontents

\section{Introduction}\label{S:Intro}

By definition, a \emph{domain} is a ring~$R$ without zero divisors, and an \emph{atomic decomposition} of an element~$a\in R$ is a product decomposition $a=a_1\cdots a_n$ where each~$a_i$ is an \emph{atom}, that is, $a_i$ is nonzero, non-unit, and is not a product of two non-units.

This paper presents two observations on atomic decompositions of elements in domains, focused on the concept of \emph{similarity \pup{of elements}} in a domain (see Section~\ref{S:Basic}).

\subsection{Order of the similarity classes of factors}\label{Su:OrdSim}
Cohn defines in \cite[Corollary~1]{Cohn1963a} a \emph{unique factorization domain} (UFD for short) as a domain in which every nonzero element~$a$ can be written as a product of atoms, and if $a=x_1\cdots x_m=y_1\cdots y_n$ are such factorizations, then $m=n$ and there exists a permutation~$\gs$ of $\set{1,\dots,n}$ such that each~$x_i$ is similar to~$y_{\gs i}$.
Notable classes of UFDs are \emph{rings with a weak algorithm} (see Cohn \cite[Corollary~2]{Cohn1963b}), and \emph{\pup{not necessarily commutative} principal ideal domains} (PIDs for short; see Jacobson~\cite[Ch.~3]{Jaco1943}).
Those classes of rings both contain the one of \emph{Euclidean domains} (see~\cite{Cohn2006a}), of which a prominent subclass consists of the polynomial rings~$\kk[\xx]$ where~$\kk$ is a division ring and~$\xx$ a central indeterminate.

In the latter case, we answer the question whether the similarity classes can be realized in any prescribed order, by constructing a skew field~$\kk$ such that~$\kk[\xx]$ fails that statement, even for factorizations of length~$2$.
The construction requires a division ring~$\kk$ with non-conjugate elements~$a$ and~$b$ for which \emph{Sylvester's equation} $at-tb=1$ has no solution in~$\kk$.
Such an example was already constructed by Berghout in~\cite{Bergh1968}, with $|a|=|b|$.
Our construction has $|a|\neq|b|$; that inequality will be required in subsequent work.

\subsection{Transferring atomic decompositions through similarity}
\label{Su:TransfSim}

Cohn states in \cite[Proposition~2.2]{Cohn1963a} that \emph{for any similar elements~$a$ and~$b$ in a domain, any factorization of~$a$ gives rise to a factorization of~$b$}.
We shall verify that this statement does not hold as a rule.

In order to pin the concept in question, let us state a definition.

\begin{definition}\label{D:SimSplit}
A domain~$R$ is \emph{similarity-splitting} if for all $a_1,a_2,b\in R$, if $a_1a_2\sim b$ then there are $b_1,b_2\in R$ such that $b=b_1b_2$ and each $a_i\sim b_i$.
\end{definition}

\begin{note}
It suffices to verify the given condition for~$a_1$, $a_2$, and~$b$ all nonzero and non-units (see Corollary~\ref{C:Rsimpp1}).
\end{note}

In view of Definition~\ref{D:SimSplit}, Cohn's \cite[Proposition~2.2]{Cohn1963a} can be restated as
\begin{equation}\label{Eq:SimSplit}
\text{\em Every domain is similarity-splitting.}
\end{equation}

We shall give a sufficient condition for Cohn's conclusion to hold for a domain~$R$.
This condition will be stated ``every right ideal with cyclic quotient is principal'' (see Proposition~\ref{P:StrCyc2Pal2Decomp}).
It holds if~$R$ is either commutative or weak B\'ezout (see Proposition~\ref{P:StrCycl}).

Nonetheless, we shall construct in Section~\ref{S:NonSimSplit} non-similarity-splitting domains.
In each of those examples, there are similar elements one of them an atom and the other not a product of atoms (Example~\ref{Ex:StrCycl}) or a product of two atoms one of them central (Example~\ref{Ex:StrCycl2}).
Thus both examples disprove Cohn's original claim \cite[Proposition~2.2]{Cohn1963a} in a strong sense.
In each of them there exists a non principal right ideal with cyclic quotient (see also Remark~\ref{Rk:NonStrCycl}).

\subsection{Some history pertaining to \cite[Proposition~2.2]{Cohn1963a}}\label{Su:Cohn2.2}
It is not clear to us whether Cohn was indeed aware of a gap in his original argument (or statement).
We could not locate any place, in the literature, where such a gap would be mentioned explicitly.

Cohn's \cite[Proposition~2.2]{Cohn1963a} is also stated, with the same argument, as Proposition~1.3.4 in Cohn's 2006 monograph~\cite{Cohn2006a}.

On the other hand, Exercise~3.2(2) in Cohn~\cite{Cohn1985} (restated without change in Cohn \cite[Exercise~3.1(3)]{Cohn2006a}) asks the reader to prove that in any $2$-fir (or, more generally, $2$-Hermite domain)~$R$, $R/I\cong R/aR$ (as right $R$-modules) implies that~$I$ is principal, whenever~$I$ is a right ideal of~$R$ and~$a\in R$.

Exercise~3.2(6) in Cohn~\cite{Cohn1985} (restated without change in Cohn \cite[Exercise~3.1(8)]{Cohn2006a}) asks the reader to prove that for any similar elements~$c$ and~$c'$ in a $2$-fir~$R$, every decomposition $c=ab$ gives rise to a decomposition~$c'=a'b'$ with~$a'$ similar to~$a$ and~$b'$ similar to~$b$; moreover, the same exercise asks ``Does this hold in more general rings?".

Further references on factorization in rings are Cohn~\cite{Cohn1973b,Cohn1973a}, Smertnig~\cite{Smer2016}, Facchini and Fassina~\cite{FacFas2018}, Bell \emph{et al.}~\cite{BBNS2023}.

\subsection{Use of AI}\label{Su:AI}

The author asked the AI system \texttt{ChatGPT} (OpenAI, web interface; conversation dated October 30, 2026) to check the proofs for correctness and suggest stylistic changes.
On November~1, 2026, further checks by more powerful versions of both \texttt{ChatGPT} and \texttt{Claude} (courtesy of a colleague with access to those models) were performed.
A few of their suggestions were adopted, mostly typographical and stylistic corrections, but also the references to Berghout~\cite{Bergh1968}, Lam and Leroy~\cite{LamLeroy2004}.
No AI assistance was involved in the discovery process leading to this paper.

\section{Basic concepts}\label{S:Basic}

All our rings will be associative and unital, but not necessarily commutative.
All ring homomorphisms and all subrings will be unital.
We will set $\ZZ^+\eqdef\set{0,1,2,\dots}$.

Elements~$a$ and~$b$ in a ring~$R$ are \emph{right similar}, written $a\simr b$, if\linebreak $R/aR\cong R/bR$ as right $R$-modules.
\emph{Left similarity} is the symmetric concept (i.e., $a\siml b$ if $R/Ra\cong R/Rb$ as left $R$-modules).
It is well known, and nontrivial (see Cohn \cite[Corollary~1]{Cohn1963a}%
\footnote{
The Lemma on page~315 in Cohn~\cite{Cohn1963a} is incorrect as stated (e.g., there are counterexamples to the stated equivalence between~$a$ and~$b$ having a unit HCLF and $aR+bR=R$; moreover, $a$ and~$b$ should be both assumed nonzero).
We leave the fix to the reader as an exercise.
},
referring to Fitting~\cite{Fitt1936}), that if~$R$ is a domain, then right and left similarity of~$a$ and~$b$ are equivalent; we shall then say that~$a$ and~$b$ are similar, and write $a\sim b$.
Associated elements (i.e., $b=uav$ for units~$u$ and~$v$) are similar, but not conversely.%
\footnote{Smertnig incorrectly claims in \cite[Example~4.2(2)]{Smer2016} that for every field~$\FF$, the generators~$x$ and~$y$ of the free associative algebra~$\FF\seq{x,y}$ are similar.
On the other hand, the elements $1+xy$ and $1+yx$ are similar but not associated.
}
For an automorphism~$\tau$ of a ring~$R$ and an indeterminate~$\xx$, the \emph{skew polynomial ring} $R[\xx;\tau]$ consists of all finite sums of ``monomials'' of the form $\xx^na$, where $a\in R$ and $n\in\ZZ^+$, subjected to the commutation rules
 \begin{equation}\label{Eq:axtaua}
 a\xx=\xx\tau(a)\quad\text{for}\quad a\in R\,.
 \end{equation}
It is well known that for every automorphism~$\tau$ of a division ring~$\kk$, the ring $\kk[\xx;\tau]$ satisfies Ore's condition on both sides (see Cohn \cite[\S~7.3]{Cohn2003b}), thus it embeds into a unique division ring of quotients~$\kk(\xx;\tau)$, whose elements can be written in each of the forms~$uv^{-1}$ and~$v^{-1}u$ where $u,v\in R[\xx;\tau]$ and $v\neq0$.
Moreover, $\kk(\xx;\tau)$ embeds in turn into the division ring~$\ls{\kk}{\xx;\tau}$ of all \emph{skew Laurent series}, consisting of all formal linear combinations of the form $\sum_{n\geq m}\xx^na_n$ where $m\in\ZZ$ and each $a_n\in\kk$, subjected to the commutation rules~\eqref{Eq:axtaua}.

In case~$\tau$ is the identity, that construction specializes to the polynomial ring~$\kk[\xx]$ with one \emph{central} indeterminate~$\xx$.
By iterating the latter construction and forming directed colimits, we get the polynomial rings $\kk[\xx_i\mid i\in I]$ over a collection\linebreak $\vecm{\xx_i}{i\in I}$ of central indeterminates.
Those rings are domains satisfying Ore's condition on both sides (see Proposition~7.3.4 and Theorem~7.3.5 in Cohn~\cite{Cohn2003b}), thus giving rise to their division rings of quotients $\kk\vecm{\xx_i}{i\in I}$, which are non-commutative generalizations of the classical rational function fields.

An \emph{absolute value} on a division ring~$\kk$ is a map $\kk\to\RR^+$, $x\mapsto|x|$ such that the statements $|x+y|\leq|x|+|y|$, $|xy|=|x|\cdot|y|$, and ($|x|=0\Leftrightarrow x=0$), hold whenever $x,y\in\kk$.
A \emph{valued division ring} is a pair $(\kk,|{}_-|)$ where~$|{}_-|$ is an absolute value on~$\kk$.
If the relation $|x+y|\leq|x|+|y|$ is replaced by $|x+y|\leq\max\set{|x|,|y|}$ we say that~$|{}_-|$ is \emph{non-Archimedean}.

Weak B\'ezout domains are introduced in Cohn~\cite{Cohn1963a}, as domains in which any two principal right ideals with nonzero intersection have an intersection and a sum which are both principal.
They are also identified as \emph{$2$-firs} in~\cite[\S~2]{Cohn1970}.
Cohn proves in~\cite{Cohn1963a} that being weak B\'ezout is a left-right symmetric condition.
By Williams~\cite{Will1968}, \emph{a domain is weak B\'ezout if{f} any two principal right ideals with non-zero intersection generate a principal right ideal}.

We set $R^{\times}\eqdef R\setminus\set{0}$, for every domain~$R$.

\section{Getting started on similarity}\label{S:Similar}

All the results in this section belong to the folklore.
We provide proofs for convenience.

At the most elementary level, right similarity can be expressed by a positive primitive statement:

\begin{lemma}\label{L:Rsimpp}
Elements~$a$ and~$b$ in a unital ring~$R$ are right similar if{f} there exist $u,v,p,q,r,s\in \nobreak R$ such that $ua=bp$, $vb=aq$, $1-vu=ar$, and $1-uv=bs$.
\end{lemma}

\begin{proof}
The homomorphisms $R/aR\to R/bR$ are exactly those of the form $x+aR\mapsto ux+bR$ where $ua\in bR$, and the homomorphisms $R/bR\to R/aR$ are exactly those of the form $y+bR\mapsto vy+aR$ where $vb\in aR$.
Such homomorphisms are mutually inverse if{f} $1-vu\in aR$ and $1-uv\in bR$.
\end{proof}

A ring is \emph{directly finite} if it satisfies the implication $xy=1\Rightarrow yx=1$.
In particular, \emph{if a ring is either finite or a domain, then it is directly finite}.

\begin{corollary}\label{C:Rsimpp1}
The following statements hold, for every element~$a$ in a directly finite ring~$R$:
\begin{enumerater}
\item\label{asim0}
$a\simr0$ if{f} $a=0$.

\item\label{asim1}
$a\simr1$ if{f}~$a$ is a unit of~$R$.

\end{enumerater}
\end{corollary}

\begin{proof}
We verify the nontrivial directions.

\eqref{asim0}.
By Lemma~\ref{L:Rsimpp} (with $b=0$), there are $u,v\in R$ such that $ua=0$ and $uv=1$.
Since~$R$ is directly finite, $vu=1$; whence $a=vua=0$.

\eqref{asim1}.
By definition, $a\simr1$ if{f} $aR=R$, that is, $a$ is a left unit.
Since~$R$ is directly finite, $a$ is a unit.
\end{proof}

\begin{corollary}\label{C:Rsimpp2}
Let $f\colon R\to S$ be a ring homomorphism and let $a,b\in R$.
If~$a$ and~$b$ are right similar within~$R$, then $f(a)$ and~$f(b)$ are right similar within~$S$.
\end{corollary}

\begin{lemma}\label{L:1-ab}
For any elements~$a$ and~$b$ in a ring~$R$, the elements $1-ab$ and $1-ba$ are \pup{both right and left} similar within~$R$.
\end{lemma}

\begin{proof}
There are mutually inverse isomorphisms 
 \begin{align*}
 R/(1-ba)R\to R/(1-ab)R\,,&&x+(1-ba)R\mapsto ax+(1-ab)R\,,\\
 R/(1-ab)R\to R/(1-ba)R\,,&&x+(1-ab)R\mapsto bx+(1-ba)R\,.
 \end{align*}
This settles right similarity.
Left similarity is symmetric.
\end{proof}

The following lemma entails Cohn's observation, stated in \cite[Exercise~3.1(15)]{Cohn2006a}, that ``\emph{In any skew polynomial ring, similar elements have the same degree}''.

\begin{lemma}\label{L:Cohn2006aExo}
For any division ring~$\kk$ and $p\in\kk[\xx]^{\times}$, the dimension of $\kk[\xx]/p\kk[\xx]$, viewed as a right $\kk$-vector space, is the degree of~$p$.
\end{lemma}

\begin{proof}
Setting $m\eqdef\deg{p}$, the finite sequence $\vec{x}\eqdef\vecm{\xx^k+p\kk[\xx]}{0\leq k<m}$ is (trivially) linearly independent over~$\kk$.
Applying right Euclidean division by~$p$ in~$\kk[\xx]$, we obtain that~$\vec{x}$ generates $\kk[\xx]/p\kk[\xx]$.
\end{proof}

\begin{corollary}\label{C:Sim2Deg}
Let~$\kk$ be a division ring, let $a\in\kk$, and let $p\in\kk[\xx]$.
Then the following are equivalent:
\begin{enumeratei}

\item\label{asimp}
$\xx-a\sim p$;

\item\label{lmconj}
$p=\gl\cdot(\xx-a')$ for some conjugate~$a'$ of~$a$ within~$\kk$;

\item\label{rmconj}
$p=(\xx-a')\cdot\gl$ for some conjugate~$a'$ of~$a$ within~$\kk$.
\end{enumeratei}
\end{corollary}

\begin{proof}
The implications \eqref{lmconj}$\Leftrightarrow$\eqref{rmconj} and \eqref{lmconj}$\Rightarrow$\eqref{asimp} are straightforward.
Now assume~\eqref{asimp} and set $R\eqdef\kk[\xx]$.
By Lemma~\ref{L:Cohn2006aExo}, the degree of~$p$ is~$1$, and both right $\kk$-vector spaces $A\eqdef R/(\xx-a)R$ and $P\eqdef R/pR$ are one-dimensional, with an isomorphism $\gf\colon A\to P$.
Let $\gl\in\kkm$ such that $\gf(1+(\xx-a)R)=\gl+pR$.
Then
 \[
 0=\gf(0)=\gf\pI{(1+(\xx-a)R)(\xx-a)}=(\gl+pR)(\xx-a)
 =\gl(\xx-a)+pR\,,
 \]
whence $\gl(\xx-a)\in pR$.
Since $\gl(\xx-a)$ and~$p$ both have degree~$1$, it follows that there exists $\xi\in\kkm$ such that $p=\gl\cdot(\xx-a)\cdot\xi$.
Therefore,
 \begin{equation*}
 p=\gl\xi(\xx-\xi^{-1}a\xi)=(\xx-\gl a\gl^{-1})\gl\xi\,.\tag*{\qed}
 \end{equation*}
\renewcommand{\qed}{}
\end{proof}

\section{An unsolvable Sylvester's equation with nonconjugate parameters}\label{S:UnsolvSylv}

As observed in Cohn \cite[Lemma~2.3]{Cohn1973c}, existence and uniqueness of the solution to Sylvester's equation $at-tb=c$ play an important role in the study of spectral properties of matrices over division rings.
In this section we shall relate the existence of a solution to Sylvester's equation $at-tb=1$, in a division ring~$\kk$, to factorizations of elements of~$\kk[\xx]$ with degree~$2$.

Our next observation%
\footnote{Proposition~\ref{P:Johnson}\eqref{Johnson} is contained in Johnson \cite[Theorem~1]{Johns1944}.
We include a (straightforward) proof for convenience.}
records a few straightforward instances of existence and uniqueness of a solution for Sylvester's equation.

\begin{proposition}\label{P:Johnson}
The following statement holds, for any elements~$a$, $b$, $c$ in a division ring~$\kk$:
\begin{enumerater}
\item\label{UniqSylv}
If~$a$ and~$b$ are not conjugate, then Sylvester's equation
 \begin{equation}\label{Eq:SylvEqn}
 at-tb=c
 \end{equation}
has at most one solution in~$\kk$.
\end{enumerater}
Moreover, each of the following conditions implies that~\eqref{Eq:SylvEqn} has a unique solution:
\begin{enumerater}
\setcounter{enumi}{1}
\item\label{Johnson}
$a$ and~$b$ are not conjugate and $\kk$ is finite-dimensional over its center~$\FF$;

\item\label{CplSyl}
$|a|\neq|b|$ and~$\kk$ is a complete valued division ring.

\end{enumerater}

\end{proposition}

\begin{proof}
\eqref{UniqSylv}.
The map $f\colon\kk\to\kk$, $x\mapsto ax-xb$ is $\FF$-linear.
If~$a$ and~$b$ are not conjugate, then~$f$ is one-to-one; the uniqueness statement follows.

\eqref{Johnson}.
By~\eqref{UniqSylv}, $f$ is one-to-one.
If~$\kk$ is finite-dimensional over~$\FF$, then~$f$ must thus be an $\FF$-linear automorphism.

\eqref{CplSyl}.
{F}rom $|a|\neq|b|$ it follows that~$a$ and~$b$ are not conjugate; thus, by~\eqref{UniqSylv}, $f$ is one-to-one.
If $|a|>|b|$, then the series $a^{-1}c+a^{-2}cb+a^{-3}cb^2+\cdots$ converges to a solution of~\eqref{Eq:SylvEqn}.
Similarly, if $|a|<|b|$ then the series $-cb^{-1}-acb^{-2}-a^2cb^{-3}-\cdots$ converges to a solution of~\eqref{Eq:SylvEqn}.
\end{proof}

A particular case of Proposition~\ref{P:Johnson} is the division ring~$\HH$ of all quaternions: \emph{for any non-conjugate quaternions~$a$ and~$b$, Sylvester's equation $at-tb=c$ has a unique solution whenever $c\in\HH$}.

Finding a division ring with non-conjugate elements~$a$ and~$b$ such that Sylvester's equation $at-tb=1$ has no solution is less straightforward.
Berghout does this in~\cite{Bergh1968}, by producing a complete valued division ring~$\kk$ of characteristic~$0$, with non-conjugate elements~$a$, $b$ such that there is no~$t$ with $at-tb=1$.

Berghout's example has~$|a|=|b|$.
We shall construct a class of examples with $|a|\neq|b|$; by Proposition~\ref{P:Johnson}, this forces~$\kk$ to be incomplete.
Our construction will be based on the following definition, which is a slightly modified version of a suggestion by George Bergman.

\begin{definition}\label{D:VdMde}
A triple of the form $(\kk,a,\tau)$, where~$\kk$ is a division ring, $a\in\nobreak\kk$, and~$\tau$ is an automorphism of~$\kk$, is a \emph{Vandermonde triple} if the row vectors
 \[
 u_n\eqdef(a^n,\tau(a^n),\tau^2(a^n),\dots)\,,\qquad\text{for }n\in\ZZ^+\,,
 \]
are right linearly independent over~$\kk$ (i.e., $\sum_nu_n\gl_n=0$ entails all $\gl_n=0$).
\end{definition}

For a family $\vecm{\gl_n}{n\in\ZZ^+}$, with finite support, of elements in~$\kk$, the condition $\sum_nu_n\gl_n=0$ means that all elements~$\tau^l(a)$, for $l\in\ZZ^+$, are left roots of the polynomial $p(\xx)\eqdef\sum_n\xx^n\gl_n$.
Now if~$\kk$ is commutative and~$p\neq0$, then~$p$ has only finitely many roots in~$\kk$.
The argument can easily be reversed.
We thus get:
 \begin{quote}\em
 If~$\kk$ is commutative, then $(\kk,a,\tau)$ is a Vandermonde triple if{f} all elements~$\tau^l(a)$, for $l\in\ZZ^+$, are pairwise distinct.
 \end{quote}
In particular, we get the following examples, suggested to the author by George Bergman.

\begin{example}\label{Ex:VdMchar0}
Let~$\FF$ be a field and set $\kk\eqdef\FF(s)$ with~$s$ an indeterminate.
The following statements hold:
\begin{itemize}
\item
Suppose that~$\FF$ has characteristic~$0$ and let~$\tau_1$ be the unique automorphism of~$\FF(s)$ over~$\FF$ such that $\tau_1(s)=s+1$.
Then $(\kk,s,\tau_1)$ is a Vandermonde triple.

\item
Suppose that~$\FF$ has a nonzero element~$q$ which is not a root of unit and let~$\tau_2$ be the unique automorphism of~$\FF(s)$ over~$\FF$ such that $\tau_2(s)=qs$.
Then $(\kk,s,\tau_2)$ is a Vandermonde triple.

\end{itemize}

\end{example}

The statement and proof of our next result benefited substantially from a mail exchange with George Bergman.

\begin{proposition}\label{P:UnsolvSylv}
Let~$(\kk,a,\tau)$ be a Vandermonde triple.
Endow the division ring $\kk(b;\tau)$, where~$b$ is an indeterminate, with its $b$-adic absolute value.
Then Sylvester's equation $at-tb=1$ has no solution in $\kk(b;\tau)$, whereas $|b|<|a|$.
\end{proposition}

\begin{proof}
Our argument, of non-existence of a solution to Sylvester's equation, will take place within the tower $\kk[b;\tau]\hookrightarrow\kk(b;\tau)\hookrightarrow\ls{\kk}{b;\tau}$.
Recall that the $b$-adic absolute value on $\ls{\kk}{b;\tau}$ is defined by setting $|0|=0$ and
 \[
 |u|\eqdef2^{-n}\ \text{where }n\text{ minimal such that }u_n\neq0\,,\quad
 \text{whenever }u\neq0\,.
 \]
In particular, $|a|=1$ whereas $|b|=1/2$.

As in the proof of Proposition~\ref{P:Johnson}\eqref{CplSyl}, we observe the following.

\begin{sclaim}
The power series $u\eqdef\sum_{n\geq0}a^{-n-1}b^n$ is the unique element of $\ls{\kk}{b;\tau}$ such that $au-ub=1$.
\end{sclaim}

Thus, in order to complete the proof of Proposition~\ref{P:UnsolvSylv}, it suffices to verify that~$u$ does not belong to~$\kk(b;\tau)$.

Suppose otherwise.
This means that $ux\in\kk[b;\tau]$ for some $x\in\kk[b;\tau]^{\times}$.
Write $x=\sum_{0\leq k\leq m}b^kx_k$, with each $x_k\in\kk$.
Then $ux=\sum_{n\geq0}b^ny_n$ where each $y_n\in\kk$.
The condition $ux\in\kk[b;\tau]$ thus means that there exists an integer~$m'$, which may be taken above~$m$, such that $y_n=0$ whenever $n\geq m'$.
Now for each $n\geq m$,
 \[
 y_n=\sum\nolimits_{0\leq l\leq m}\tau^n(a^{l-n-1})x_l\,,
 \]
thus for all $n\geq m'$, we get, multiplying the equation $y_n=0$ by $\tau^n(a^{n+1})$ to the left,
 \begin{equation}\label{Eq:nullyn}
 \sum\nolimits_{0\leq l\leq m}\tau^n(a^l)x_l=0\,.
 \end{equation}
Setting $x'_l\eqdef\tau^{-m'}(x_l)$ whenever $0\leq l\leq m$, and applying the automorphism~$\tau^{-m'}$ to both sides of~\eqref{Eq:nullyn}, we get
 \begin{equation}\label{Eq:lfloornullyn}
 \sum\nolimits_{0\leq l\leq m}\tau^n(a^l)x'_l=0\quad\text{whenever }n\in\ZZ^+\,.
 \end{equation}
Now we can write~\eqref{Eq:lfloornullyn} in the form
 \[
 \sum\nolimits_{0\leq l\leq m}(a^l,\tau(a^l),\tau^2(a^l),\dots)x'_l=
 (0,0,0,\dots)\,,
 \] 
which, since $(\kk,a,\tau)$ is a Vandermonde triple, implies that each $x'_l=0$, thus $x_l=0$; in contradiction with our assumption that $x\neq0$.
\end{proof}

The following observation provides a large supply of Vandermonde triples.

\begin{proposition}\label{P:ConstrVDM}
For any division ring~$D$, let $\kk\eqdef D\vecm{s_i}{i\in\ZZ}$ \pup{the~$s_i$ are central indeterminates}, and let~$\tau$ be the unique automorphism of~$\kk$ over~$D$ mapping each~$s_i$ to~$s_{i+1}$.
Then $(\kk,s_0,\tau)$ is a Vandermonde triple.
\end{proposition}

\begin{proof}
We need to verify that the vectors $u_k\eqdef(s_0^k,s_1^k,s_2^k,\dots)$, for $k\in\ZZ^+$, are right linearly independent over~$\kk$.
This means that for any sequence $\vecm{x_n}{n\in\ZZ^+}$ with finite support, of elements in~$\kk$, the conjunction of all equations
 \begin{equation}\label{Eq:DsiiZZ}
 x_0+s_nx_1+s_n^2x_2+\cdots=0\,,
 \end{equation}
for $n\in\ZZ^+$, implies that all $x_n=0$.
Since $D[s_i\mid i\in\ZZ]$ satisfies Ore's condition, we may clear a common denominator from the~$x_i$ in~\eqref{Eq:DsiiZZ} and thus assume that all~$x_i$ are polynomials.
Now pick $n\in\ZZ^+$ such that~$s_n$ is absent from any of the~$x_i$.
Clearly, \eqref{Eq:DsiiZZ} implies that all $x_i=0$.
\end{proof}

Bringing together Propositions~\ref{P:UnsolvSylv} and~\ref{P:ConstrVDM}, we get

\begin{corollary}\label{C:ConstrVDM}
Every division ring can be embedded into a non-Archimedean valued division ring~$\kk$, with elements~$a$ and~$b$ of distinct absolute value, such that Sylvester's equation $at-tb=1$ has no solution in~$\kk$.
\end{corollary}

\begin{note}
In Corollary~\ref{C:ConstrVDM} both orderings of $(|a|,|b|)$ are possible: if $|b|<|a|$ (as it is the case in the construction of Proposition~\ref{P:UnsolvSylv}), replace the triple $(\kk,a,b)$ by $(\kk^{\op},b,a)$, where~$\kk^{\op}$ denotes the opposite division ring of~$\kk$, endowed with the same absolute value as~$\kk$.
The new values of~$a$ and~$b$ satisfy $|a|<|b|$.
\end{note}

Lam and Leroy already observed in \cite[Lemma~10.1]{LamLeroy2004} the correspondence between existence of solutions to certain Sylvester's equations and decomposition properties of polynomials of degree~$2$ over skew fields.
Our next statement is of similar nature.

\begin{lemma}\label{L:abnotba}
Let~$a$, $b$, $c$, $d$ be elements in a division ring~$\kk$ with $b\neq c$.
If
 \begin{equation}\label{Eq:x-bacd}
 (\xx-b)(\xx-a)=(\xx-c)(\xx-d)\qquad\text{within }\kk[\xx]\,,
 \end{equation}
then the equation $at-tb=1$ has a solution in~$\kk$.
\end{lemma}

\begin{proof}
Expanding~\eqref{Eq:x-bacd} yields $a+b=c+d$ and $ba=cd$.
Since $b\neq c$, the element $s\eqdef c-b=a-d$ is nonzero, whence
 \[
 ba=cd=(b+s)(a-s)=ba-bs+sa-s^2\,,
 \]
which simplifies to $sa-bs=s^2$.
Multiplying both sides of that equation by~$s^{-1}$ on the left and on the right, it follows that $as^{-1}-s^{-1}b=1$ which, writing
$t\eqdef s^{-1}$, gives the desired relation.
\end{proof}

\begin{proposition}\label{P:abnotba}
Let~$a$ and~$b$ be nonconjugate elements in a division ring~$\kk$ such that the equation $at-tb=1$ has no solution in~$\kk$.
Then there are no polynomials $p,q\in\kk[\xx]$ such that
 \begin{equation}\label{Eq:condsonp,q}
 p\sim\xx-a\,,\quad q\sim\xx-b\,,\quad\text{and}\quad (\xx-b)(\xx-a)=pq\,.
 \end{equation}
\end{proposition}

\begin{proof}
Let $p,q\in\kk[\xx]$ satisfy~\eqref{Eq:condsonp,q}.
By Corollary~\ref{C:Sim2Deg}, there are $\xi,\eta\in\kkm$, a conjugate~$c$ of~$a$, and a conjugate~$d$ of~$b$ within~$\kk$ such that $p=(\xx-c)\cdot\xi$ and $q=\eta\cdot(\xx-d)$.
It follows that $(\xx-b)(\xx-a)=(\xx-c)\xi\eta(\xx-d)$.
Equating the coefficients of~$\xx^2$ yields $\xi\eta=1$, whence $(\xx-b)(\xx-a)=(\xx-c)(\xx-d)$.
Since~$a$ is conjugate to~$c$ but not to~$b$, we get $b\neq c$.
By Lemma~\ref{L:abnotba}, we get a contradiction.
\end{proof}

\section{Similarity-splitting domains}
\label{S:StrCycl}

The (faulty, see  Section~\ref{S:NonSimSplit}) proof of Cohn \cite[Proposition~2.2]{Cohn1963a} argues there that a chain of cyclic submodules from~$R$ to~$aR$, gives rise, \emph{via} the isomorphism $R/aR\cong R/bR$, to a chain from~$R$ to~$bR$, in which corresponding factors are isomorphic (this is the argument of the proof of Proposition~\ref{P:StrCyc2Pal2Decomp}).
As Cohn pointed out in~\cite{Cohn1965}, the ``$R$ to~$aR$'' above should be followed by ``with strictly cyclic quotients'', where a module~$M$ is strictly cyclic if $M\cong\nobreak R/cR$ for some $c\in R^{\times}$.

\begin{definition}\label{D:CyclWit}
For a right ideal~$I$ of a ring~$R$, a \emph{cyclicity witness for~$I$} is a triple $(a,b,c)$ of elements in~$R$ such that
 \begin{equation}\label{Eq:CyclWit}
 1-ab\in I\qquad\text{and}\qquad(ax\in I\Leftrightarrow x\in cR
 \quad\text{whenever }x\in R)\,.
 \end{equation}
\end{definition}

The following lemma is a direct translation of our definitions; we omit the proof.

\begin{lemma}\label{L:CyclWit}
For a right ideal~$I$ of~$R$, the module~$R/I$ is cyclic if{f} there exists a cyclicity witness for~$I$.
\end{lemma}

\begin{lemma}\label{L:CyclWitn}
Let~$(a,b,c)$ be a cyclicity witness for a right ideal~$I$ of a ring~$R$.
Then $1-ba\in cR$ and $I=(1-ab)R+acR$.
\end{lemma}

\begin{proof}
{F}rom $1-ab\in I$ it follows that $a(1-ba)=(1-ab)a\in I$, thus, since $(a,b,c)$ is a cyclicity witness, $1-ba\in cR$.

The right ideal $I'\eqdef(1-ab)R+acR$ is obviously contained in~$I$.
Let, conversely, $x\in I$.
Since $(1-ab)x\in I'\subseteq I$ and $x\in I$ we get $abx\in I$, thus, by~\eqref{Eq:CyclWit}, $bx\in cR$; whence $abx\in acR\subseteq I'$.
Since $(1-ab)x\in I'$, it follows that\linebreak $x=abx+(1-ab)x\in I'$.
\end{proof}

Being a cyclicity witness can be written as a positive primitive statement:

\begin{lemma}\label{L:CyclWitExist}
A triple $(a,b,c)$ of elements in a ring~$R$ is a cyclicity witness for some right ideal of~$R$ if{f} $1-ba\in cR$.
\end{lemma}

\begin{proof}
By Lemma~\ref{L:CyclWitn}, every cyclicity witness $(a,b,c)$ satisfies $1-ba\in cR$.
Let, conversely, $a,b,c,d\in R$ such that $ba+cd=1$.
We shall verify that $(a,b,c)$ is a cyclicity witness for the right ideal $I\eqdef(1-ab)R+acR$.
The conditions $1-ab\in I$ and $ac\in I$ are trivially satisfied.
It remains to verify that every element~$x$ of~$R$ such that $ax\in I$ belongs to~$cR$.
There are $y,z\in R$ such that $ax=(1-ab)y+acz$.
In particular, $y=as$ for $s\eqdef x+by-cz$, whence, setting $t\eqdef ds+z$,
 \[
 ax=(1-ab)as+acz=a(1-ba)s+acz=acds+acz=act\,.
 \]
It follows that $ba(x-ct)=0$, that is, $(1-cd)(x-ct)=0$; whence $x\in cR$.
\end{proof}

Bringing together Lemmas~\ref{L:CyclWitn} and~\ref{L:CyclWitExist}, we obtain

\begin{corollary}\label{C:CyclWitExist}
For any ring~$R$, the following conditions are equivalent:
\begin{enumeratei}
\item\label{allISC}
every right ideal of~$R$ with cyclic quotient is principal;

\item\label{bacdSC}
for all $a,b,c,d\in R$, $ba+cd=1$ implies that the right ideal $(1-ab)R+acR$ is principal.

\end{enumeratei}
\end{corollary}

We do not know whether Condition~\eqref{bacdSC} above is left-right symmetric.

\section{Fixing Cohn's original statement}\label{S:FixingCohn}

This section will be centered on a sufficient condition, for a ring~$R$, to be similarity-splitting (see Definition~\ref{D:SimSplit}); that is, to satisfy the conclusion of \cite[Proposition~2.2]{Cohn1963a}.

Our next result cannot be considered as new: the commutative case is trivial, and the weak B\'ezout case follows from Cohn \cite[Exercise~3.1(3)]{Cohn2006a}.
Since the proof is short we include it for convenience.

\begin{proposition}\label{P:StrCycl}
Let~$I$ be a right ideal in a domain~$R$, with $R/I$ cyclic.
If~$R$ is either commutative or weak B\'ezout, then~$I$ is principal.
\end{proposition}

\begin{note}
For every commutative ring~$R$, every ideal~$I$ is the annihilator of any generator of~$R/I$.
Hence, $R/I\cong R/J$ if{f} $I=J$.
In particular, if~$R$ is a domain, then $a\sim b$ if{f} $a$ and~$b$ are associated.
\end{note}

\begin{proof}
Throughout the proof fix a cyclicity witness $(a,b,c)$ for~$I$.
By Lemma~\ref{L:CyclWitn}, $I=(1-ab)R+acR$.

Suppose first that~$R$ is commutative.
Then $c=(1-ab)c+acb\in I$ thus $cR\subseteq I$.
Conversely, for each $x\in I$, $ax\in I$ thus $x\in cR$; whence $I\subseteq cR$.
Therefore, $I=cR$ is principal.

Suppose now that~$R$ is weak B\'ezout.
Pick $d\in R$ such that $1=ba+cd$ (see Lemma~\ref{L:CyclWitn}).
If $ba=1$, then, since~$R$ is a domain, $ab=1$, thus $I=acR$ so we are done.
If $a=0$, then $I=R$ so we are done again.
Suppose from now on that $ba\neq1$, that is, $cd\neq0$, and $a\neq0$.
Since $0\neq acd=a(1-ba)=(1-ab)a$, we get $acR\cap(1-ab)R\neq\set{0}$, thus, by our assumption, $\set{ac,1-ab}$ generates a principal right ideal, which must be~$I$.
\end{proof}

In particular, the assumptions of Proposition~\ref{P:StrCycl} are satisfied if~$R$ is a ring with weak algorithm.
This holds, in particular, if~$R$ is a free associative algebra over a field (see Cohn~\cite{Cohn1963b}).

Domains in which every right ideal with cyclic quotient is principal yield a corrected version of \cite[Proposition~2.2]{Cohn1963a}, with Cohn's original argument as proof:

\begin{proposition}\label{P:StrCyc2Pal2Decomp}
Let~$R$ be a domain.
If every right ideal of~$R$ with cyclic quotient is principal, then $R$ is similarity-splitting.
\end{proposition}

\begin{proof}
Let $a_1,a_2,b\in R$ such that $a_1a_2\sim b$.
We need to find $b_1,b_2\in R$ such that $b=b_1b_2$ and each $a_i\sim b_i$.
By Corollary~\ref{C:Rsimpp1}, we may assume that~$a_1\neq0$. 

Set $a\eqdef a_1a_2$, and pick mutually inverse isomorphisms $\ga\colon R/bR\to R/aR$ and $\gb\colon R/aR\to R/bR$.
Then $\gb[a_1R/aR]=I/bR$ for a unique right ideal~$I$ of~$R$ containing~$bR$, and further,
 \[
 R/a_1R\cong(R/aR)/(a_1R/aR)\cong(R/bR)/\gb[a_1R/aR]
 =(R/bR)/(I/bR)\cong R/I\,,
 \]
so, by our assumption on~$R$, $I$ is principal; so $I=b_1R$ for some $b_1\in R$ with $a_1\sim b_1$.
Since $a_1\neq0$ and by Corollary~\ref{C:Rsimpp1}, $b_1\neq0$.
{F}rom $bR\subseteq b_1R$ it follows that $b=b_1b_2$ for some $b_2\in R$.
Now
 \[
 R/b_2R\cong b_1R/b_1b_2R=I/bR\cong\ga[I/bR]=
 a_1R/a_1a_2R\cong R/a_2R\,,
 \]
so $a_2\sim b_2$.
\end{proof}

Note that being similarity-splitting is a left-right symmetric condition, while the property ``every right ideal with cyclic quotient is principal'' does not look left-right symmetric \emph{a priori}.

\begin{remark}[No addition of similarity classes]\label{Rk:simNotAdd}
Let~$R$ be a domain.
In view of Proposition~\ref{P:StrCyc2Pal2Decomp}, it is tempting to endow the set $\cS=\setm{[x]}{x\in R}$, of similarity classes of all elements of~$R$, with a binary operation~$+$, meant as an ``addition'', by setting
 \begin{equation}\label{Eq:AddSimCl}
 [a]+[b]=[ab]\,,\quad\text{whenever }a,b\in R\,.
 \end{equation}
However, if this were possible, then the statement
 \begin{equation}\label{Eq:AddSim}
 y\sim z\Rightarrow xy\sim xz\,,\quad\text{whenever }x,y,z\in R
 \end{equation}
would hold.
We shall now verify that this is not the case, by proving the following statement:
 \begin{sclaim}\label{Cl:xyxynotsim}
 Denote by~$\FF_2$ the two-element field and by~$x$, $y$ non-commuting
 indeterminates.
 Then $x(yx+1)\not\sim x(xy+1)$ within $\FF_2\seq{x,y}$\,,
 \end{sclaim}

\begin{scproof}
Suppose that $x(yx+1)\sim x(xy+1)$ within $\FF_2\seq{x,y}$.
Let~$R$ be the ring of all lower triangular matrices over~$\FF_2$, with elements 
 \[
 a\eqdef\begin{bmatrix}0 & 0\\ 0 & 1\end{bmatrix}\,,\quad
 b\eqdef\begin{bmatrix}0 & 0\\ 1 & 1\end{bmatrix}
 \]
and let $f\colon\FF_2\seq{x,y}\to R$ be the unique homomorphism carrying $(x,y)$ to $(a,b)$.
By Corollary~\ref{C:Rsimpp2}, $f\pI{x(yx+1)}\simr f\pI{x(xy+1)}$ within~$R$; that is,
 \[
 a(ba+1)\simr a(ab+1)\quad\text{within }R\,.
 \]
Now $a(ba+1)=0$ whereas $a(ab+1)=a+b\neq0$.
Since~$R$ is finite, it is directly finite, which, by Corollary~\ref{C:Rsimpp1}, leads to a contradiction.
\end{scproof}
Now we can complete our argument: since $xy+1\sim yx+1$ (see Lemma~\ref{L:1-ab}), \eqref{Eq:AddSim} fails for the domain~$\FF_2\seq{x,y}$.
Therefore, \eqref{Eq:AddSimCl} also fails for that domain.

\end{remark}

\section{Construction of non-similarity-splitting domains}\label{S:NonSimSplit}

Our first example already disproves Cohn's \cite[Proposition~2.2]{Cohn1963a} in a strong sense.

\begin{examplepf}\label{Ex:StrCycl}
For any field~$\FF$, we construct an $\FF$-algebra~$Z$ without zero divisors, with similar elements in~$Z$ such that one is an atom whereas the other is not a product of atoms.
In particular, $Z$ is not similarity-splitting.
\end{examplepf}

\begin{proof}
The polynomial ring $\FF[s]$, where~$s$ is an indeterminate, is a subring of the rational function field~$\FF(s)$.

Endow $V\eqdef\FF(s)[\xx]$ with its natural structure of $\FF(s)$-vector space and denote by~$\End{V}$ its endomorphism ring.
We define endomorphisms~$\ol{a}$ and~$\ol{b}$ of~$V$ by the rules $\ol{b}(1)=0$, $\ol{b}(\xx^{n+1})=\xx^n$, and $\ol{a}(\xx^n)=\xx^{n+1}$ whenever $n\in\ZZ^+$.
Observe that $\ol{b}\circ\ol{a}=\id_V$ whereas $(\ol{a}\circ\ol{b})(u)=u-u(0)$ whenever $u\in V$.%
\footnote{
It is stated in Lam \cite[Exercise~11.9]{Lam2001} (credited there to Jacobson) that the assignment $x\mapsto\ol{a}$, $y\mapsto\ol{b}$ gives rise to an isomorphism from $\FF\seq{x,y}/{(1-yx)}$ onto the subalgebra $\FF\seq{\ol{a},\ol{b}}$ of~$\End{V}$ generated by $\{\ol{a},\ol{b}\}$.}

Denote by~$\FF(s)\seq{a,b}$ the free associative $\FF(s)$-algebra on non-commuting indeterminates~$a$ and~$b$, by ($\FF(s)\seq{a,b}\twoheadrightarrow\FF(s)$, $z\mapsto z(1,1)$) the evaluation map, and by $\gf\colon\FF(s)\seq{a,b}\twoheadrightarrow\FF(s)\seq{\ol{a},\ol{b}}$ the unique $\FF(s)$-algebra homomorphism carrying $(a,b)$ to $(\ol{a},\ol{b})$.

Denote by~$Q$ the set of all $f\in\End{V}$ for which there are $\xi\in\FF(s)[\xx,\xx^{-1}]$ and a nonnegative integer~$m$ such that
 \[
 f(u)=\xi u\quad\text{whenever}\quad u\in\xx^m\FF(s)[\xx]\,.
 \]
It is then easy to verify that the Laurent polynomial~$\ell(f)\eqdef\xi$ is uniquely determined by~$f$.
Moreover, $Q$ is a subring of~$\End{V}$ containing $\set{\ol{a},\ol{b}}$ (thus the range of~$\gf$) and~$\ell$ is an $\FF(s)$-algebra homomorphism from~$Q$ to~$\FF(s)[\xx,\xx^{-1}]$.
By composing~$\ell$ with the evaluation morphism $\eps\colon\FF(s)[\xx,\xx^{-1}]\to\FF(s)$, $\xi\mapsto\xi(1)$, we obtain an $\FF(s)$-algebra homomorphism $\pi=\eps\ell\colon Q\to\FF(s)$, $f\mapsto\ell(f)(1)$.

Denote by~$T$ the subring of~$Q$ consisting of the endomorphisms that carry~$\FF[s,\xx]=\FF[s][\xx]$ into itself.
Alternatively, $T$ consists of the elements of~$Q$ that can be represented, in the basis $\vecm{\xx^k}{k\in\ZZ^+}$, by column-finite matrices with entries in~$\FF[s]$.

Our example will be $Z\eqdef\gf^{-1}[T]$.
It contains the subset $\set{a,b,c,s}$ of~$\FF(s)\seq{a,b}$ where $c\eqdef s^{-1}(1-ba)$.

\setcounter{claim}{0}

\begin{claim}\label{Cl:ell(T)}
For each $f\in T$, $\ell(f)\in\FF[s,\xx,\xx^{-1}]$ and $\pi(f)\in\FF[s]$.
\end{claim}

\begin{cproof}
Let $m\in\ZZ^+$ such that $f(u)=\ell(f)u$ whenever $u\in\xx^m\FF(s)[\xx]$.
Taking $u=\xx^m$ yields $\ell(f)=\xx^{-m}f(\xx^m)$.
From $f\in T$ it follows that $f(\xx^m)\in\FF[s,\xx]$.
Hence, $\ell(f)\in\xx^{-m}\FF[s,\xx]\subseteq\FF[s,\xx,\xx^{-1}]$, and hence the element $\pi(f)=\nobreak\ell(f)(1)$ belongs to~$\FF[s]$.
\end{cproof}

The various structures involved, with $\FF$-algebra homomorphisms between them, are represented in Figure~\ref{Fig:StrCycl} (unlabeled hook arrows $\hookrightarrow$ represent inclusion maps).
Ignoring the two inclusion arrows $\FF(s)\hookrightarrow\FF(s)\seq{a,b}$ and $\FF[s]\hookrightarrow Z$, the diagram commutes.
Moreover, $\pi\circ\gf\res_{\FF(s)}=\id_{\FF(s)}$.
The bottom right arrow $\pi\res_T\colon T\to\FF[s]$ is justified by Claim~\ref{Cl:ell(T)}.

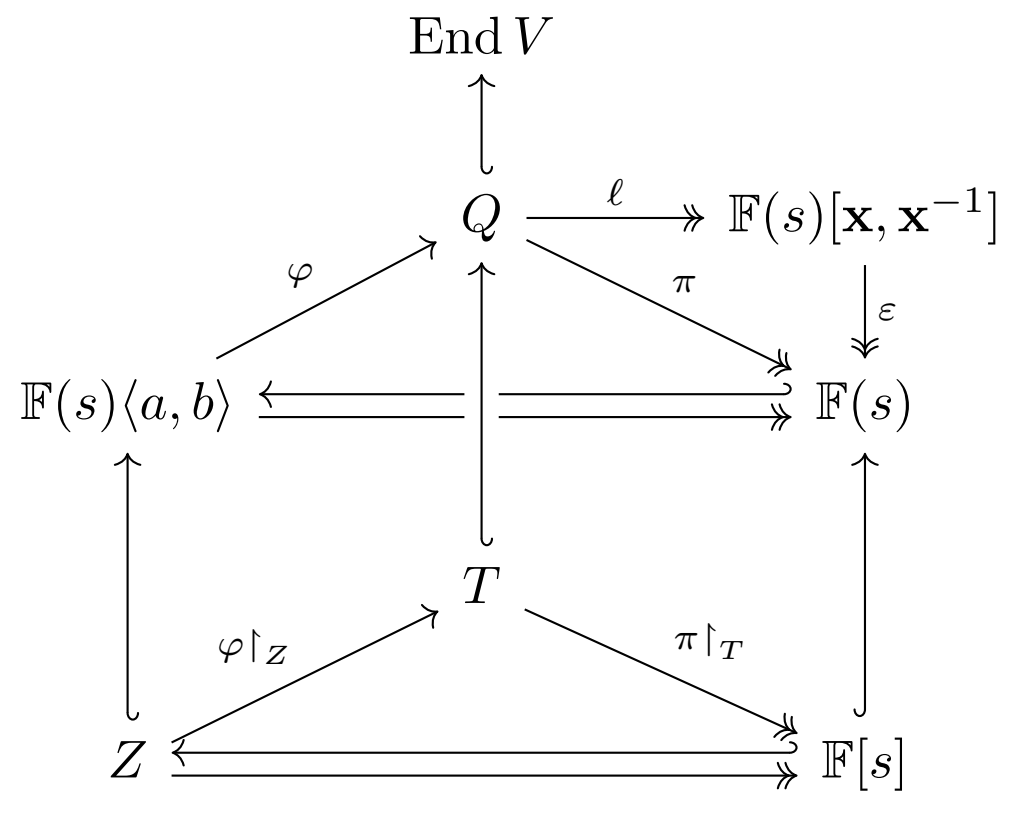
\begin{figure}[htb]
\begin{tikzcd}
\centering
& \End{V} & \\
& Q \arrow[r,"\ell",two heads]\arrow[rd,"\pi",two heads]\arrow[u,hook'] & \FF(s)[\xx,\xx^{-1}]\arrow[d,two heads,"\eps"]\\
\FF(s)\seq{a,b}\arrow[ru,"\gf"]\arrow[rr,two heads,shift right=.5ex] &  & \FF(s)\arrow[ll,hook',shift right=.5ex]\\
& T\arrow[uu,hook',crossing over]\arrow[rd,"\pi\res_T",two heads] \\
Z\arrow[uu,hook']\arrow[ru,"\gf\res_Z"]\arrow[rr,two heads,shift right=.5ex] & &
\FF[s]\arrow[ll,hook',shift right=.5ex]
\arrow[uu,hook]
\end{tikzcd}
\caption{Illustrating the construction of
Example~\ref{Ex:StrCycl}}
\label{Fig:StrCycl}
\end{figure}

\begin{claim}\label{Cl:RandZZ}
The relation $\pi\gf(z)=z(1,1)$ holds whenever $z\in\FF(s)\seq{a,b}$.
Moreover, $z\in Z$ implies that $z(1,1)\in\FF[s]$.
In particular, $Z\cap\FF(s)=\FF[s]$.
\end{claim}

\begin{cproof}
From $\ell(\ol{a})=\xx$ and $\ell(\ol{b})=\xx^{-1}$ it follows that
 \[
 \pi\gf(a)=\pi(\ol{a})=\ell(\ol{a})(1)=\xx(1)=1=a(1,1)\,.
 \]
Likewise, $\ell(\ol{b})=\xx^{-1}$ and $\pi\gf(b)=1=b(1,1)$.
Further, $\pi\gf(1)=1$.
This completes the proof of the first part of our claim.

For all $z\in Z$, $\gf(z)$ belongs to~$T$, thus, using Claim~\ref{Cl:ell(T)} and the paragraph above, $z(1,1)=\pi\gf(z)\in\FF[s]$.
The final statement of Claim~\ref{Cl:RandZZ} follows.
\end{cproof}

%

\begin{claim}\label{Cl:1/2(1-ab)}
Set $c_{\gl}\eqdef\gl(1-ba)$ and $d_{\gl}\eqdef\gl(1-ab)$ whenever $\gl\in\FF(s)$.
Then $c_{\gl}\in Z$ for all~$\gl$, whereas $d_{\gl}\in Z$ if{f} $\gl\in\FF[s]$.
\end{claim}

\begin{cproof}
{F}rom $\ol{b}\ol{a}=\id_V$ it follows that $c_{\gl}\in\ker\gf\subseteq Z$.
Since $\gf(d_{\gl})(1)=\gl 1$, $d_{\gl}\in Z$ implies that $\gl\in\FF[s]$.
The converse is trivial.
\end{cproof}

\begin{claim}\label{Cl:1-abindec}
$1-ab$ is an atom of~$Z$.
\end{claim}

\begin{cproof}
Let $u,v\in Z$ such that $1-ab=uv$.
Since $1-ab$ is an atom of~$\FF(s)\seq{a,b}$, either~$u$ or~$v$ belongs to~$\FF(s)$, thus $1-ab=uv=vu$ and we may assume that $u\in\FF(s)$; whence $u\in\FF[s]$ by Claim~\ref{Cl:RandZZ}.
Since $v=u^{-1}(1-ab)$ and by Claim~\ref{Cl:1/2(1-ab)}, it follows that $u^{-1}\in\FF[s]$, so~$u$ is a unit of~$\FF[s]$; that is, $u\in\FFm$.
\end{cproof}

\begin{claim}\label{Cl:1-baverydec}
$1-ba$ is not a product of atoms in~$Z$.
\end{claim}

\begin{cproof}
Since $1-ba$ is an atom of~$\FF(s)\seq{a,b}$, any factor of $1-ba$ is either a scalar (in~$\FF(s)$) or~$c_{\gl}$ for some $\gl\in\FF(s)^{\times}$ (see Claim~\ref{Cl:1/2(1-ab)}).
Thus any product decomposition of $1-ba$ must include at least one~$c_{\gl}$.
Now $c_{\gl}=s\cdot c_{s^{-1}\gl}$ is not an atom of~$Z$.
\end{cproof}

The desired conclusion follows from Lemma~\ref{L:1-ab} together with Claims~\ref{Cl:1-abindec} and~\ref{Cl:1-baverydec}.
\end{proof}

Our next example is a variant of Example~\ref{Ex:StrCycl}, now with the second element a product of two atoms.

\begin{examplepf}\label{Ex:StrCycl2}
For any field~$\FF$, we construct an $\FF$-algebra~$S$ without zero divisors, with similar elements in~$S$ such that one is an atom whereas the other is a product of two atoms one of which is central.
In particular, $S$ is not similarity-splitting.
\end{examplepf}

\begin{proof}
Our example will be constructed as an $\FF$-subalgebra of the domain~$Z$ constructed in Example~\ref{Ex:StrCycl}.
We shall keep the notation of the proof of that example.

Let~$t$ be a new central indeterminate and set
 \[
 \gS\eqdef\bigoplus\vecm{\FF s^mt^n}{(m,n)\in\ZZ\times\ZZ\,,\ m+n\geq0}\,,
 \]
a proper $\FF$-subalgebra of $\FF[s,s^{-1},t,t^{-1}]$, thus, \emph{a fortiori}, of~$\FF(s)[t,t^{-1}]$, containing~$\FF[s]$.
Our example will be
 \begin{equation}\label{Eq:defnS}
 S\eqdef\setm{z\in Z}{z(1+t,1)\in\gS}\,,
 \end{equation}
which is an $\FF$-subalgebra of~$Z$ (thus of $\FF(s)\seq{a,b}$).
Note that in~\eqref{Eq:defnS}, $z$ belongs to~$Z$, thus to $\FF(s)\seq{a,b}$, and thus $z(1+t,1)$ belongs to $\FF(s)[t,t^{-1}]$.
Trivially, $\set{1,a,b}\subseteq S$.
Reusing the elements $c_{\gl}\eqdef\gl(1-ba)$ and $d_{\gl}\eqdef\gl(1-ab)$ from the proof of Claim~\ref{Cl:1/2(1-ab)} in Example~\ref{Ex:StrCycl}, we state:

\setcounter{claim}{0}

\begin{claim}\label{Cl:SinterswFF}
For every $\gl\in\FF(s)$, the following statements hold:
\begin{enumerater}
\item\label{glinS}
$\gl\in S$ if{f} $\gl\in\FF[s]$;

\item\label{gltinS}
$\gl t\in\gS$ if{f} $\gl\in s^{-1}\FF[s]$;

\item\label{cglinS}
$c_{\gl}\in S$ if{f} $\gl\in s^{-1}\FF[s]$;

\item\label{dglinS}
$d_{\gl}\in S$ if{f} $\gl\in\FF[s]$.

\end{enumerater}
\end{claim}

\begin{cproof}
\eqref{glinS}.
We verify the nontrivial direction.
By our definition~\eqref{Eq:defnS}, $\gl\in S$ if{f} $\gl\in\gS$.
If that holds, then we can write
 \begin{equation}\label{Eq:glingS}
 \gl=\sum\nolimits_{i<N}\xi_is^{m_i}t^{n_i}\,,\ \text{with }N\in\ZZ^+
 \text{ least possible},\text{ each }
 m_i+n_i\geq0\text{ and }\xi_i\in\FFm\,.
 \end{equation}
Since~$\gl\in\FF(s)$, it follows, arguing within $\FF(s)[t,t^{-1}]$, that all $n_i=0$.
Therefore, all $m_i\geq0$ and $\gl=\sum_{i<N}\xi_is^{m_i}\in\FF[s]$.

\eqref{gltinS}.
We verify the nontrivial direction.
Since $\gl t\in\gS$ we can write
 \begin{equation}\label{Eq:gltingS}
 \gl t=\sum\nolimits_{i<N}\xi_is^{m_i}t^{n_i}\,,\ \text{with }N\in\ZZ^+
 \text{ least possible},\text{ each }
 m_i+n_i\geq0\text{ and }\xi_i\in\FFm\,.
 \end{equation}
Since $\gl\in\FF(s)$, it follows, again arguing within $\FF(s)[t,t^{-1}]$, that all $n_i=1$; whence all $m_i\geq-1$.
Furthermore, evaluating~\eqref{Eq:gltingS} at $(s,1)$ yields $\gl=\sum_{i<N}\xi_is^{m_i}$.
Since all $m_i\geq-1$, we get $\gl\in s^{-1}\FF[s]$.

\eqref{cglinS}.
We have seen in Claim~\ref{Cl:1/2(1-ab)} of the proof of Example~\ref{Ex:StrCycl} that $c_{\gl}\in Z$ whenever $\gl\in\FF(s)$.
Hence, by~\eqref{Eq:defnS}, $c_{\gl}\in S$ if{f} $c_{\gl}(1+t,1)\in\gS$ if{f} $\gl t\in\gS$.
By~\eqref{gltinS}, this is equivalent to $\gl\in s^{-1}\FF[s]$.

\eqref{dglinS}.
We have seen in Claim~\ref{Cl:1/2(1-ab)} of the proof of Example~\ref{Ex:StrCycl} that $d_{\gl}\in Z$ if{f} $\gl\in\FF[s]$.
Hence, by~\eqref{Eq:defnS}, $d_{\gl}\in S$ if{f} $\gl\in\FF[s]$ and $\gl t\in\gS$.
Apply~\eqref{gltinS} and observe that $\FF[s]\subset s^{-1}\FF[s]$.
\end{cproof}

\begin{claim}\label{Cl:1-abbairred}
The element $1-ab$ is an atom of~$S$, whereas $1-ba$ is the product of two atoms, one of them central.
\end{claim}

\begin{cproof}
Since $1-ab$ is an atom of~$\FF(s)\seq{a,b}$, any decomposition of that element as a product of two factors within~$S$ has the form $1-ab=\gl^{-1}\cdot\gl(1-ab)$ where $\gl\in\FF(s)^{\times}$, $\gl^{-1}\in S$, and $\gl(1-ab)\in S$.
By Claim~\ref{Cl:SinterswFF}, $\gl$ and~$\gl^{-1}$ should both belong to~$\FF[s]$, thus $\gl\in\FFm$ is a unit of~$S$.
This completes the proof of the atomicity of~$1-ab$ within~$S$.

As to the second statement of Claim~\ref{Cl:1-abbairred} and since $1-ba=s\cdot s^{-1}(1-ba)$, with~$s$ central and $\set{s,s^{-1}(1-ba)}\subseteq S$ (see Claim~\ref{Cl:SinterswFF}) it suffices to verify that~$s$ and $s^{-1}(1-\nobreak ba)$ are both atoms of~$S$.
The first statement follows immediately from Claim~\ref{Cl:SinterswFF}\eqref{glinS} together with the atomicity of~$s$ within~$\FF[s]$.
Since $s^{-1}(1-ba)$ is an atom of~$\FF(s)\seq{a,b}$, any decomposition of that element as a product of two factors within~$S$ has the form $s^{-1}(1-ba)=\gl\cdot\gl^{-1}s^{-1}(1-ba)$ where $\gl\in\FF[s]^{\times}$ (because of Claim~\ref{Cl:SinterswFF}\eqref{glinS}) and $\gl^{-1}s^{-1}\in s^{-1}\FF[s]$ (because of Claim~\ref{Cl:SinterswFF}\eqref{cglinS}).
The latter condition means that $\gl^{-1}\in\FF[s]$.
Therefore, $\gl$ is a unit of~$\FF[s]$, that is, an element of~$\FFm$.
\end{cproof}

Since $1-ab\sim s\cdot s^{-1}(1-ba)$ within~$S$ (see Lemma~\ref{L:1-ab}) with $1-ab$, $s$ (which is also central), and $s^{-1}(1-ba)$ all atoms of~$S$, the proof is complete.
\end{proof}

\begin{remark}\label{Rk:NonStrCycl}
By Proposition~\ref{P:StrCyc2Pal2Decomp}, each of the underlying rings in Examples~\ref{Ex:StrCycl} and~\ref{Ex:StrCycl2} has at least a non-principal right ideal with cyclic quotient.
Let us prove this directly, by verifying that in both constructions, the right ideal generated by $\set{1-ab,ac}$ is non-principal with a cyclic quotient.

Let us do this for the ring~$Z$ constructed in Example~\ref{Ex:StrCycl}; the proof for~$S$ (in Example~\ref{Ex:StrCycl2}) is, \emph{mutatis mutandis}, identical.
First, since $ba+sc=1$ with $\set{a,b,c,s}\subseteq Z$, it follows from Lemmas~\ref{L:CyclWitn} and~\ref{L:CyclWitExist} that $(a,b,c)$ is a cyclicity witness for the right ideal $I\eqdef(1-ab)Z+acZ$.
It thus remains to check that~$I$ is not principal.

Suppose otherwise.
There are $d,u,v,p,q\in Z$ such that
 \begin{equation}\label{Eq:uvpq0}
 d=(1-ab)(u+s^{-1}av)\,,\ 1-ab=dp\,,\text{ and }s^{-1}(1-ab)a=dq\,.
 \end{equation}
It follows that $1-ab=(1-ab)(u+s^{-1}av)p$ and $s^{-1}(1-ab)a=(1-ab)(u+s^{-1}av)q$, thus, since $\FF(s)\seq{a,b}$ is a domain and $1-ab\neq0$,
 \begin{equation}\label{Eq:uvpq}
 1=(u+s^{-1}av)p\text{ and }s^{-1}a=(u+s^{-1}av)q\,.
 \end{equation}
Since $\set{p,u+s^{-1}av}\subseteq\FF(s)\seq{a,b}$, the first relation of~\eqref{Eq:uvpq} entails~$p\in\FF(s)^{\times}$.
By the first equations in~\eqref{Eq:uvpq0}, \eqref{Eq:uvpq}, we get $d=p^{-1}(1-ab)$.
By Claim~\ref{Cl:1/2(1-ab)} of the proof of Example~\ref{Ex:StrCycl} and since $p^{-1}\in\FF(s)$, we get $p^{-1}\in\FF[s]$.
{F}rom~\eqref{Eq:uvpq} we get $p^{-1}q=s^{-1}a$, thus $q=ps^{-1}a$.
Since $ps^{-1}\in\FF(s)$, $q\in Z$, and by Claim~\ref{Cl:RandZZ} of the proof of Example~\ref{Ex:StrCycl}, we get $ps^{-1}\in\FF[s]$; whence (since $p^{-1}\in\FF[s]$) $s^{-1}\in p^{-1}\FF[s]\subseteq\FF[s]$, a contradiction.
\end{remark}

\section{Acknowledgment}\label{S:Acknow}
The author thanks George Bergman for having generously shared many comments and insights on the early versions of this paper, which led to many improvements, especially in Section~\ref{S:UnsolvSylv}.


\providecommand{\noopsort}[1]{}\def\cprime{$'$}
  \def\polhk#1{\setbox0=\hbox{#1}{\ooalign{\hidewidth
  \lower1.5ex\hbox{`}\hidewidth\crcr\unhbox0}}} \def\cprime{$'$}
  \def\cprime{$'$} \def\cprime{$'$} \def\cprime{$'$} \def\cprime{$'$}
\providecommand{\bysame}{\leavevmode\hbox to3em{\hrulefill}\thinspace}
\providecommand{\MR}{\relax\ifhmode\unskip\space\fi MR }
\providecommand{\MRhref}[2]{%
  \href{http://www.ams.org/mathscinet-getitem?mr=#1}{#2}
}
\providecommand{\href}[2]{#2}

\end{document}